\documentclass[reqno]{amsart}

\usepackage{upref}
\usepackage{amsfonts, amsmath}
\newtheorem{theorem}{Theorem}

\numberwithin{equation}{section}
\keywords{Semilinear elliptic equations, Fixed point theorem, Krasnoselskii}
\thanks{This research was supported by CAPES, BRAZIL}

\begin{document}
\title[Semilinear elliptic equations and fixed points]{ 
Semilinear elliptic equations and fixed points\\
}

\author[Cleon S. Barroso]{Cleon S. Barroso}

\address{
   Departamento de Matematica,
   Universidade Federal do Cear\'a,\newline 
   Campus do Pici, Bl.914, Fortaleza-Ce, 60455-760, Brazil}
\email{cleonbar@mat.ufc.br}   

\date{}

\maketitle 

\allowdisplaybreaks

\begin{abstract}
In this paper, we deal with a class of semilinear elliptic equation in a bounded domain $\Omega\subset\mathbb{R}\sp N$, $N\geq 3$, with $C\sp{1,1}$ boundary.
Using a new fixed point result of the Krasnoselskii's type for the sum of two operators, an existence principle of strong solutions is proved. We give two examples where the nonlinearity can be critical.
\end{abstract}

\section{Introduction}
Consider the boundary value problem for a semilinear elliptic equation:
\begin{equation*}\label{eqn:1}\left \{
\begin{split}
&-\Delta u+\lambda u=f(x,u,\mu)\quad\mbox{ in } \Omega, \\
& u=0\quad \mbox{  on } \partial\Omega,
\end{split}\right.
\tag{$P\sb{\lambda,\mu}$}
\end{equation*}
where $\Omega$ is a bounded domain in $\mathbb{R}\sp N$ $(N\geq 3)$ with $C\sp{1,1}$ boundary, $\Delta$ is the standard Laplace operator, $\lambda$ is a parameter close to 0, and $f:\Omega\times\mathbb{R}\times\mathbb{R}\sp{+}\to\mathbb{R}$ is a Caratheodory function. Elliptic equations such as (\ref{eqn:1}) have been studied extensively for several years (see \cite{1,2,3,10} and the references therein). According to the behaviour of $f(\cdot,u)$, topological methods may come to be more appropriate to the solvability of (\ref{eqn:1}) than variational tools. Specially when $f(\cdot,u)$ possesses critical growth due to the lack of compactness of Sobolev embedding.
\par
In this paper, we study the problem (\ref{eqn:1}) by using the fixed point methods. A new version of Krasnoselskii's fixed point theorem for the sum of two operators given for the author \cite{6} is used in order to establish an existence principle for this problem. Here we will be interested in the existence of strong solution to (\ref{eqn:1}), that is, a function $u\in W\sp{2,2}(\Omega)\cap W\sp{1,2}\sb 0 (\Omega)$ such that $-\Delta u(x)+\lambda u(x)=f(x,u(x))$ for almost every $x\in\Omega$ and fulfilling the boundary condition.
\par
A special case of (\ref{eqn:1}) is the problem
\begin{equation*}\label{eqn:2}\left \{
\begin{split}
&-\Delta u+\lambda u=\mu u\sp p+h(x)\quad\mbox{ in }\Omega,\\
&u=0\quad\mbox{ on } \partial\Omega,
\end{split}\right.
\tag{$Q\sb{\lambda,\mu}$}
\end{equation*}
where $p>1$, $\mu>0$ and $h\in L\sp 2(\Omega)$. When $\lambda=0$ and $\mu=1$ this problem relates to an open question (cf. \cite[pg. 124]{1}): For any $h$ in $L\sp 2(\Omega)$, does $(\ref{eqn:2})$ have infinitely many solutions? Bahri \cite{8} has given a partial answer. He proved for $p\in ]2,2\sp *[$, that there is an open dense set of $h$ in $W\sp{-1,2}(\Omega)$ for which $(\ref{eqn:2})$ possesses an infinite number of distinct weak solutions. In \cite{9} the authors showed that there exists $1<p\sb N<(N+2)/(N-2)$ such that for any $h\in L\sp 2(\Omega)$ and $p\in (1,p\sb N)$ the previous equation possesses infinitely many distinct solutions.
Here, we show for any $h\in L\sp 2(\Omega)$, $h\not\equiv 0$, and $p>1$ if $N=3$, or $1<p\leq N/(N-4)$ if $N>4$, that there exists a positive constant $\mu\sb{p,h}$ sucht that $(\ref{eqn:2})$ possesses a nontrivial strong solution if $\mu\in (0,\mu\sb{p,h}]$. Beyond this, if $h\in L\sp 2(\Omega)$ satisfies $h\geq 0$ almost everywhere in $\Omega$, then the solution is positive. Notice that if $N>4$ then $(N+2)/(N-2)< N/(N-4)$. Thus, the technique used in \cite{8,9} cannot be applied to solving $(\ref{eqn:2})$ when $p=2\sp *$ and $N=3$, or $p\sb N\leq p\leq N/(N-4)$ if $N>4$.
\par
The paper is organized in four sections, including the Introduction. The main result and some preliminaries are presented in Section 2. In Section 3 we prove the main result. In the last Section we study the problem $(\ref{eqn:2})$ and another special case of $(\ref{eqn:1})$ concerning to an eigenvalue problem.

\section{Main Result}
Throughout this paper $\Omega$ denotes a bounded domain in $\mathbb{R}\sp N$ $(N\geq 3)$ with $C\sp{1,1}$ boundary. In this section, let $E$ denotes the following sapce $W\sp{2,2}(\Omega)\cap W\sp{1,2}\sb 0 (\Omega)$, where $W\sp{2,2}(\Omega)$ is the usual Sobolev space and $W\sp{1,2}\sb 0 (\Omega)$ is the closure of 
$C\sb 0\sp{\infty}(\Omega)$ in the norm $\|u\|\sb{1,2,\Omega}$. It is well-known that the mapping $\phi:u\mapsto -\Delta u$ is one-to-one from $E$ onto $L\sp 2(\Omega)$. Moreover, there is a constant $C>0$ such that 
\begin{equation}\label{eqn:3}
\|u\|\sb{2,2}\leq C\|\Delta u\|\sb 2,
\end{equation}
for all $u\in E$, see \cite[Lemma 9.17]{7}. Write now $\lambda\sp *=\|\phi\sp{-1}\|\sb{\mathcal{L}\sp 2(\Omega)}\sp{-1}$.
\par
Given a Caratheodory function $f:\Omega\times\mathbb{R}\times\mathbb{R}\sp{+}\to\mathbb{R}$, the assumptions below will be posited in what follows.
\begin{itemize}
\item[$(\mathcal{H}\sb 1)$] there is $\mu>0$ such that 
$N\sb{f\sb{\mu}}(\phi\sp{-1})(B\sb R)\subset B\sb R$, for some $R>0$ ;
\item[$(\mathcal{H}\sb 2)$] there is $\mu>0$ such that 
$N\sb{f\sb{\mu}}(\phi\sp{-1})\sp{+}(B\sb R)\subset B\sb R$, for some $R>0$;
\item[$(\mathcal{H}\sb 3)$] $f\sb{\mu}(x,u)=f(x,u,\mu)\geq 0$ if $u\geq 0$;
\item[$(\mathcal{H}\sb 4)$] $p>1$ if $N=3$, or $1<p<N/(N-4)$ if $N>4$,
\end{itemize}
where $N\sb{f\sb{\mu}}(u)=f(\cdot,u,\mu)$ (resp. $N\sb{f\sb{\mu}} (u)\sp{+}=f\sb{\mu}(\cdot,u\sp{+})$) and $B\sb R =B\sb R (0)\subset L\sp 2(\Omega)$.
\par
In this paper we prove the following theorem.

\begin{theorem}\label{trm:1} Let $0\leq\lambda\leq\lambda\sp *$. Suppose  $f:\Omega\times\mathbb{R}\times\mathbb{R}\sp{+}\to\mathbb{R}$ is a Caratheodory function. Then the following holds:
\begin{itemize}
\item[$(a)$] $(\ref{eqn:1})$ has a strong solution if $(\mathcal{H}\sb 1)$ is fulfilled.
\item[$(b)$] $(\ref{eqn:1})$ has a nonnegative strong solution if $(\mathcal{H}\sb 2)$-$(\mathcal{H}\sb 3)$ are verified.
\end{itemize} 
\end{theorem}

This improves upon the recent work in \cite{6}, Theorem $4.1$, which instead of $(\mathcal{H}\sb 1)$ requires the hypothesis $N\sb f\phi\sp{-1}(B\sb{RC})\subset B\sb R$ with $C>0$ satisfying (\ref{eqn:3}) and some $R>0$ . As we will see in Section $4$, in some concrete situations the assumption $(\mathcal{H}\sb 1)$, or $(\mathcal{H}\sb 2)$, seems to be more efficient than this hypothesis.

\par
The proof of theorem \ref{trm:1} is based on the following version asymptotic of the Krasnoselskii's fixed point theorem for the sum of two operators given for the author in \cite{6}. 

\begin{theorem}\label{trm:2} Suppose $X$ reflexive and assume that $B\in\mathcal{L}(X)$ with $\|B^p\|\leq 1$, $p\geq 1$, is a dissipative operator on $X$. If $A:X\to X$ is a weakly continuous mapping such that $A(B_R)\subseteq B_R$, for some $R>0$, then there is $y\in B_R$ such that $Ay+By=y$. 
\end{theorem}

\section{proof of theorem \ref{trm:1}}
The proof is similar to the Theorem $4.1$ of \cite{6} and we include here for the sake of completeness. Let $B$ be the linear operator in $L\sp 2(\Omega)$ given by $Bu=-\lambda\phi^{-1}(u)$ with $0<\lambda\leq \lambda\sp *$. Next, define $A:B\sb R\to L^2(\Omega)$ by
\begin{equation*}
A(u)=\left\{
\begin{split}
 &N\sb{f\sb{\mu}}(\phi^{-1}(u))\quad &\mbox{ for the case }& (a),\\
 &N\sb{f\sb{\mu}}(\phi^{-1}(u)\sp{+})\quad &\mbox{ for the case }& (b).
\end{split}
\right.
\end{equation*}
\paragraph{\it Claim I} $B$ is a dissipative operator in $L^2(\Omega)$ with $\|B\|\leq 1$.\\

\paragraph{\it Proof of Claim I} 
Indeed, by the choice of $\lambda$ we have $\|B\|\leq 1$. Now, because of Green's first identity we have  $(\phi(u),u)_{L^2(\Omega)}=-\int_{\Omega}\Delta u\cdot udx=\int_{\Omega}|\nabla u|^2dx\geq 0$, for all $u\in E$. Using this we see that $(Bu,u)_{L^2(\Omega)}\leq 0$, for all $u\in L^2(\Omega)$, which implies that $B$ is a dissipative operator in $L\sp 2(\Omega)$. This proves the claim I. 
\\
\paragraph{\it Claim II} $A$ is weakly continuous in $B\sb R$.\\

\paragraph{\it Proof of Claim II} It is enough to show the claim for the item $(b)$, see \cite{6}.
Clearly $A$ is well defined and by $(\mathcal{H}\sb 2)$ maps $B_R$ into itself. Let now $\{u_n\}$ be a sequence in $B\sb R$ such that $u_n\rightharpoonup u$ in $B\sb R$. Then, we have $\phi^{-1}(u_{n_j})\rightharpoonup \phi^{-1}(u)$ in $W^{2,2}(\Omega)$ (see, \cite[Lemma 9.17, pg. 242]{7}). Since the embedding $W^{2,2}(\Omega)\hookrightarrow L^2(\Omega)$ is compact, $\phi^{-1}(u_{n_j})\to \phi^{-1}(u)$ in $L^2(\Omega)$. Now, the Vainberg theorem \cite{4} gives $\phi\sp{-1}(u\sb n)\sp{+}\to \phi\sp{-1}(u)\sp{+}$  and $N\sb{ f\sb{\mu}}(\phi\sp{-1}(u\sb n)\sp{+})\to N\sb{f\sb{\mu}}(\phi\sp{-1}(u)\sp{+})$ in $L^2(\Omega)$. Hence 
$N\sb{f\sb{\mu}}(\phi\sp{-1}(u\sb n)\sp{+})\rightharpoonup N\sb{ f\sb{\mu}}(\phi\sp{-1}(u)\sp{+})$ in $L^2(\Omega)$. This shows that $A$ is weakly sequentially continuous in $B\sb R$. On the other hand, since $L\sp 2(\Omega)$ is reflexive and separable $(B\sb R,w)$ is a bounded and metrizable topological space, where $w$ denotes the weak topology of $L^2(\Omega)$. Owing to theorem $A6 (b)$ of \cite[pg. 370]{5} we conclude that $A$ is weakly continuous in $B\sb R$. This ends the proof of the claim.
\par
Hence, Theorem \ref{trm:2} with $X=L\sp 2(\Omega)$ and $A,B$ as above guarantees a fixed point of $A+B$ in $B\sb R$, say $u$. By using the assumption $(\mathcal{H}\sb 3)$ we may apply the weak maximum principle as stated for instance in \cite{1} to conclude that $v=\phi^{-1}(u)$ is a strong solution to $(\ref{eqn:1})$, with $v(x)\geq 0$ for almost every $x\in\Omega$. 
The $\lambda =0$ case follows from {\it Claim II} together with Schauder fixed point theorem for locally convex spaces. This establishes the Theorem \ref{trm:1}.\quad $\square$

\section{Special cases}
Theorem \ref{trm:1} has a variety of simple and practical special cases. As a first example, it can be used to solve the problem $(\ref{eqn:2})$ mentioned in the introduction.

\begin{theorem}\label{trm:3} Assume $(\mathcal{H}\sb 4)$ and $\lambda\geq 0$ as above. Then, for any $h\in L\sp 2(\Omega)$, $h\not\equiv 0$ there is $\mu\sb{p,h}>0$ such that the problem $(\ref{eqn:2})$ has a nontrivial strong solution if $\mu\in (0,\mu\sb{p,h}]$. Moreover, if $h\geq 0$ almost everywhere in $\Omega$ then the solution is positive .
\end{theorem}

\paragraph{\it Proof.} In fact, thanks to \cite[(7.30) pg. 158]{7} we see with the condition $(\mathcal{H}\sb 4)$ that
\begin{equation}\label{eqn:4}
\|\phi\sp{-1}(u)\|\sb{2p}\leq\gamma\|u\|\sb 2,
\end{equation}
for all $u\in L\sp 2(\Omega)$ and some $\gamma=\gamma\sb p(\Omega)>0$. Let $h$ be as above and define
$\mu\sb{p,h}=\frac{1}{(2\sp p\gamma\|h\|\sp{p-1}\sb 2)}$.
Next, fixing $\mu\in (0,\mu\sb{p,h}]$ and setting $\beta=(2\mu\gamma)\sp{1/(p-1)}$ we define 
$$
f(x,u,\mu)=\frac{u\sp p}{2\gamma}+\beta h(x).
$$
We are now going to verify the condition $(\mathcal{H}\sb 1)$. From $(\ref{eqn:4})$ we have
\begin{equation}\label{eqn:5}
\|N\sb{f\sb{\mu}}(\phi\sp{-1}(u))\|\sb 2\leq \frac{1}{2}\|u\|\sp p\sb 2+\beta\|h\|\sb 2,
\end{equation}
for all $u\in L\sp 2(\Omega)$. Since $\beta=(2\mu\gamma)\sp{1/(p-1)}\leq 1/2$, by $(\ref{eqn:5})$ it is immediate to verify that the condition $(\mathcal{H}\sb 1)$ holds for $R=1$. Then, Theorem \ref{trm:1} (a) yields a strong solution to $(\ref{eqn:1})$, say $u$. At last, it is easy to check that $v=(u/\beta)$ solves $(\ref{eqn:2})$. In a similar way, the estimate $(\ref{eqn:5})$ remains valid for $\|N\sb{f\sb{\mu}}(\phi\sp{-1}(u)\sp{+})\|\sb 2$ and hence the condition $(\mathcal{H}\sb 3)$ holds for $R=1$ too. So that, if $h(x)\geq 0$ for almost every $x\in\Omega$ we can apply Theorem \ref{trm:1} (b) to conclude that $u\geq 0$. Using again the weak maximum principle it follows that either $u\equiv 0$ or $u>0$. Since $h\not\equiv 0$ we get $u>0$. This completes the proof.\quad $\square$

\par 
The second application concerns to the following eigenvalue problem
\begin{equation}\label{eqn:6}\left \{
\begin{split}
&-\Delta u +\lambda u=\mu g(x,u)\quad \mbox{ in } \Omega,\\
&u=0\quad \mbox{ on } \partial \Omega,
\end{split}\right.
\end{equation}
where $g:\Omega\times\mathbb{R}\to \mathbb{R}$ is a Caratheodory function satisfying
$$
|g(x,u)|\leq a|u|\sp p + b(x),
$$
where $a>0$, $b\in L\sp 2(\Omega)$ and $p>1$.
\begin{theorem}\label{trm:4} Assume $(\mathcal{H}\sb 4)$ and $0\leq\lambda\leq\lambda\sp *$. Then, there is $\mu\sp *$ such that $(\ref{eqn:6})$ has at least one strong solution if $\mu\in (0,\mu\sp *]$.
\end{theorem}

\paragraph{\it Proof} Indeed, define $f(x,u,\mu)=\mu g(x,u)$, $\mu>0$. Then, in virtue $(\ref{eqn:4})$ one has
\begin{equation}\label{eqn:7}
\|N\sb{f\sb{\mu}}(\phi\sp{-1}(u))\|\sb 2\leq \mu (a\gamma\|u\|\sp p\sb 2+\|b\|\sb 2),
\end{equation}
for all $u\in L\sp 2(\Omega)$. Thus, if $u\in B\sb 1\subset L\sp 2(\Omega)$ and $\mu\sp *=\frac{1}{(a\gamma+\|b\|\sb 2)}$ then from $(\ref{eqn:7})$ the condition $(\mathcal{H}\sb 1)$ one follows for any $\mu\in (0,\mu\sp *]$. Appealing again to Theorem \ref{trm:1} (a), we find a strong solution to $(\ref{eqn:6})$.

\end{document}